\documentclass[10pt]{amsart}
\usepackage{palatino}
\usepackage[all]{xy,xypic}

\usepackage{amsfonts,amsmath,oldgerm,amssymb,amscd}

\newcommand{\ra}{\rightarrow}		

\newcommand{\by}[1]{\stackrel{#1}{\ra}}

\newcommand{\surj}{\ra\!\!\!\ra}	
\newcommand{\inj}{\hookrightarrow}
\newcommand{\ol}{\overline}		

\newcommand{\iso}{\by \sim}

\newtheorem{theorem}{Theorem}[section]
\newtheorem{proposition}[theorem]{Proposition}
\newtheorem{lemma}[theorem]{Lemma}
\newtheorem{definition}[theorem]{Definition}
\newtheorem{corollary}[theorem]{Corollary}
\newtheorem{question}[theorem]{Question}

	\newcommand{\gl}{\lambda}

\newcommand{\mm}{\mbox{$\mathfrak m$}}

\newcommand{\ot}{\mbox{\,$\otimes$\,}}	
\newcommand{\op}{\mbox{$\oplus$}}

\newcommand{\bp}{\begin{proposition}}
\newcommand{\ep}{\end{proposition}}
\newcommand{\bl}{\begin{lemma}}
\newcommand{\el}{\end{lemma}}
\newcommand{\bt}{\begin{theorem}}
\newcommand{\et}{\end{theorem}}
\newcommand{\bc}{\begin{corollary}}
\newcommand{\ec}{\end{corollary}}
\newcommand{\bd}{\begin{definition}}
\newcommand{\ed}{\end{definition}}

\def\rmk{\refstepcounter{theorem}\paragraph{{\bf Remark} \thetheorem}}

\def\rmk{\refstepcounter{theorem}\paragraph{{\bf Remark} \thetheorem}}

\def\proof{\paragraph{Proof}}

\def\definition{\refstepcounter{theorem}\paragraph{{\bf Definition} \thetheorem}}
\title{Existence of unimodular elements in a projective module}
\author{Manoj  K. Keshari} \address{Department of Mathematics, Indian Institute of
  Technology Bombay, Powai, Mumbai 400076, India.}
\email{keshari@math.iitb.ac.in} \author{Md. Ali Zinna}
\address{Department of Mathematics, Indian Institute of Technology
  Bombay, Powai, Mumbai 400076, India.}
\email{zinna2012@gmail.com}
\date{\today}

\begin{document}
\maketitle
\begin{abstract}
Let $R$ be an affine algebra over an algebraically closed field of
characteristic $0$ with dim$(R)=n$.  Let $P$ be a projective
$A=R[T_1,\cdots,T_k]$-module of rank $n$ with determinant $L$. Suppose
$I$ is an ideal of $A$ of height $n$ such that there are two
surjections $\alpha:P\surj I$ and $\phi:L\oplus A^{n-1} \surj I$.
Assume that either (a) $k=1$ and $n\geq 3$ or (b) $k$ is arbitrary but
$n\geq 4$ is even. Then $P$ has a unimodular element (see \ref{A1},
\ref{A2}).
\end{abstract}

\section{Introduction}
Let $R$ be a commutative Noetherian ring of dimension $n$. A classical
result of Serre \cite{se} asserts that if $P$ is a projective
$R$-module of rank $> n$, then $P$ has a unimodular element. However,
as is shown by the example of projective module corresponding to the
tangent bundle of an even dimensional real sphere, this result is best
possible in general. Therefore, it is natural to ask under what
conditions $P$ has a unimodular element when $rank(P)=n$. In
\cite{ra}, Raja Sridharan asked the following question.

\begin{question}\label{Q1}
 Let $R$ be a ring of dimension $n$ and $P$ be a projective $R$-module
 of rank $n$ with trivial determinant.  Suppose there is a surjection
 $\alpha:P\surj I$, where $I\subset R$ is an ideal of height $n$ such
 that $I$ is generated by $n$ elements. Does $P$ has a unimodular
 element?
\end{question} 

Raja Sridharan proved that the answer to Question \ref{Q1} is
affirmative in certain cases (see \cite[Theorems, 3, 5]{ra}) and
``negative'' in general.

 Plumstead \cite{p} generalized Serre's result and proved that if $P$
 is a projective $R[T]$-module of rank $> n$, then $P$ has a
 unimodular element.  Bhatwadekar and Roy \cite{br2} extended
 Plumstead's result and proved that projective
 $R[T_1,\ldots,T_r]$-modules of rank $>n$ have a unimodular element.
 Mandal \cite{ma} proved analogue of Plumstead that projective
 $R[T,T^{-1}]$-modules of rank $>n$ have a unimodular element.  In
 another direction, Bhatwadekar and Roy \cite{br1} proved that
 projective modules over $D=R[T_1,T_2]/(T_1T_2)$ of rank $>n$ have a
 unimodular element. Later Wiemers \cite{wi} extended this result and
 proved that if $D=R[T_1,\ldots,T_r]/\mathcal{I}$ is a discrete Hodge
 algebra over $R$ (here $\mathcal{I}$ is a monomial ideal), then
 projective $D$-modules of rank $>n$ have a unimodular element.

In view of results mentioned above, it is natural to ask the following
question.  Let $A$ be either a polynomial ring over $R$ or a Laurent
polynomial ring over $R$ or a discrete Hodge algebra over $R$. Let $P$
be a projective $A$-module of rank $n$. Under what conditions $P$ has
a unimodular element? We will mention two such results.

Bhatwadekar and Raja Sridharan \cite[Theorem 3.4]{brs2} proved: {\it
  Let $R$ be a ring of dimension $n$ containing an infinite field. Let
  $P$ be a projective $R[T]$-module of rank $n$. Assume $P_f$ has a
  unimodular element for some monic polynomial $f\in R[T]$. Then $P$
  has a unimodular element.}

Das and Raja Sridharan \cite[Theorem 3.4]{drs} proved:
 {\it Let $R$ be a ring of {\bf even} dimension $n$ containing
   $\mathbb{Q}$. Let $P$ be a projective $R[T]$-module of rank $n$
   with trivial determinant.  Suppose there is a surjection
   $\alpha:P\surj I$, where $I$ is an ideal of $R[T]$ of height $n$
   such that $I$ is generated by $n$ elements. Assume further that
   $P/TP$ has a unimodular element. Then $P$ has a unimodular
   element.}

Note that when $n$ is odd, the above result is not known.  Further,
the requirement in the hypothesis that $P/TP$ has a unimodular element
is indeed necessary, in view of negative answer of Question \ref{Q1}.
Motivated by Bhatwadekar-Sridharan and Das-Sridharan, we prove the following
results.

\bt (see \ref{MI}) Let $R$ be a ring of dimension $n$ containing an
infinite field. Let $P$ be a projective $A=R[T_1,\cdots,T_k]$-module of rank
   $n$. Assume $P_{f(T_k)}$ has a unimodular
  element for some monic polynomial $f(T_k)\in A$. Then $P$ has a
  unimodular element. \et
 
\bt(see \ref{many var})
Let $R$ be a ring  of {\bf even} dimension $n$ containing $\mathbb{Q}$.
Let $P$ be a projective $A=R[T_1,\cdots,T_k]$-module of rank
   $n$ with determinant $L$.  Suppose there is a surjection
   $\alpha:P\surj I$, where $I$ is an ideal of $A$ of height $n$ such
   that $I$ is a surjective image of $L\op A^{n-1}$. Further assume that
   $P/(T_1,\cdots,T_k)P$ has a unimodular element. Then $P$ has a
   unimodular element.
\et

\bt\label{laurent0} (see \ref{laurent}) 
Let $R$ be a ring of {\bf even} dimension $n$
containing $\mathbb{Q}$. Assume that height of the Jacobson radical of
$R$ is $\geq 2$. Let $P$ be a projective $R[T,T^{-1}]$-module of rank
$n$ with trivial determinant. Suppose there is a surjection
$\alpha:P\surj I$, where $I$ is an ideal of $R[T,T^{-1}]$ of height $n$
such that $I$ is generated by $n$ elements.
Then $P$ has a unimodular
element.
\et

\bt (see \ref{discrete})
 Let $R$ be a ring of {\bf even} dimension $n$ containing
$\mathbb{Q}$ and $P$ be a projective
 $D=R[T_1,T_2]/(T_1T_2)$-module of rank $n$ with determinant $L$.  Suppose there is a surjection
   $\alpha:P\surj I$, where $I$ is an ideal of $D$ of height $n$ such
   that $I$ is a surjective image of $L\op D^{n-1}$. Further, assume that
   $P/(T_1,T_2)P$ has a unimodular element. Then $P$ has a
   unimodular element.
 \et

In view of above results, we end this section with the following question.
\begin{question}
Let $R$ be a ring of {\bf even} dimension $n$ containing $\mathbb{Q}$.
\begin{enumerate}
\item Let $P$ be a projective $A=R[T,T^{-1}]$-module of rank $n$ with
trivial determinant. Suppose there is a surjection $\alpha:P\surj I$,
where $I\subset A$ is an ideal of height $n$ such that $I$
is generated by $n$ elements.  Assume further that $P/(T-1)P$ has a
unimodular element.  Does $P$ has a unimodular element?

\item Let $D=R[T_1,\ldots,T_k]/\mathcal{I}$ be a discrete Hodge algebra over $R$
and $P$ be a projective $D$-module of rank $n$ with determinant $L$.
Suppose there is a surjection $\alpha:P\surj I$, where $I$ is an ideal
of $D$ of height $n$ such that $I$ is a surjective image of $L\op
D^{n-1}$. Further, assume that $P/(T_1,\ldots,T_k)P$ has a unimodular
element. Does $P$ has a unimodular element?

\item Generalize question (1) replacing $A$ by
$R[X_1,\ldots,X_r,Y_1^{\pm 1},\ldots, Y_s^{\pm 1}]$ and $P/(T-1)P$ by
$P/(X_1,\ldots,X_r,Y_1-1,\ldots,Y_s-1)P$.
\end{enumerate}
\end{question}


 \section{Preliminaries}

\noindent{\bf Assumptions.} Throughout this paper, rings are assumed
to be commutative Noetherian and projective modules are finitely
generated and of constant rank. For a ring $A$, $\text{dim}(A)$ and
$\mathcal{J}(A)$ will denote the Krull dimension of $A$ and the
Jacobson radical of $A$ respectively. In this section we state some
results for later use.
\smallskip

\definition Let $R$ be a ring and $P$ be a projective $R$-module. An
element $p\in P$ is called \emph{unimodular} if there is a surjective
$R$-linear map $\phi:P \surj R$ such that $\phi(p)=1$.  In particular,
a row $(a_1,\cdots,a_n)\in R^n$ is unimodular (of length $n$) if there
exist $b_1,\cdots,b_n$ in $R$ such that $a_1b_1+\cdots+a_nb_n=1$. We
write $Um(P)$ for the set of unimodular elements of $P$.
\medskip


\bt \cite[Corollary 3.1]{kz}\label{uni} Let $R$ be a ring of dimension $n$,
$A=R[X_1,\cdots,X_m]$ a polynomial ring over $R$ and $P$ be a
projective $A[T]$-module of rank $n$.  Assume that $P/TP$ and $P_f$
both contain a unimodular element for some monic polynomial $f(T)\in
A[T]$.  Then $P$ has a unimodular element.  \et

The following result is a consequence of a result of Ravi Rao
\cite[Corollary 2.5]{rao} and Quillen's local-global principle
\cite[Theorem 1]{q}.
 
\bp\label{lg} Let $R$ be a ring of dimension $n$.  Suppose $n!$ is
invertible in $R$. Then all stably free $R[T]$-modules of rank $n$ are
extended from $R$.  \ep


\bl\cite[Lemma 3.3]{drs}\label{gene} Let $R$ be a ring and
$I=(a_1,\cdots,a_n)$ be an ideal of $R$, where $n$ is even.  Let $u,
v\in R$ be such that $uv=1$ modulo $I$. Assume further that the
unimodular row $(v,a_1,\cdots,a_n)$ is completable.  Then there exists
$\sigma \in M_n(R)$ with $det(\sigma)=u$ modulo $I$ such that if
$(a_1,\cdots,a_n)\sigma=(b_1,\cdots,b_n)$, then $b_1,\cdots,b_n$
generate $I$.  \el

\bp\cite[Proposition 4.2]{drs}\label{stably} Let $R$ be a ring
containing $\mathbb{Q}$ of dimension $n$ with
$\text{ht}(\mathcal{J}(R))\geq 1$.  Then any stably free $R(T)$-module
of rank $n$ is free, where $A(T)$ is obtained from $A[T]$ by inverting
all monic polynomials in $T$.  \ep


\bl\cite[Lemma 2.2]{b2} \label{bhat1}
Let $R$ be a ring and $I\subset R$ be an ideal of height $r$. Let $P$ and $Q$ be two 
projective $R/I$-modules of rank $r$ and let $\alpha:P\surj I/I^2$ and $\beta:Q\surj I/I^2$
 be surjections. Let $\psi:P\longrightarrow Q$ be a homomorphism such that ${\beta} \circ{\psi}=\alpha$. Then $\psi$ is an isomorphism.
\el

\definition\label{special} Let $R$ be a ring and $A=R[T,T^{-1}]$. We say
$f(T)\in R[T]$ is special monic if $f(T)$ is a monic polynomial
with $f(0)=1$. Write $\mathcal{A}$ for the ring obtained
from $A$ by inverting all special monic polynomials of $R[T]$.  Then it
is easy to see that $\text{dim}(\mathcal{A})=\text{dim}(R)$.
\medskip


\bt\cite[4.4, 4.6, 4.8, 4.9]{k}\label{global}
Let $R$ be a ring containing $\mathbb{Q}$ of dimension $n\geq3$ and $A=R[T,T^{-1}]$. Assume that
$\text{ht}(\mathcal{J}(R))\geq 2$. Let
$I\subset A$ be an ideal of height $n$ such that $I/I^2$ is generated by $n$ elements and
$\omega_I:(A/I)^n\surj I/I^2$ be a local orientation of $I$. Let $P$ be a projective $A$-module of rank $n$  and
$\chi:A\iso\wedge^n(P)$ an isomorphism. Then following holds.
\begin{enumerate}
 \item Suppose that the image of $(I,\omega_I)$ is zero in $E(A)$. Then $\omega_I$ is a global orientation of $I$, i.e., 
 $\omega_I$ can be lifted to a surjection from $A^n$ to $I$. 
\item  Suppose that $e(P,\chi)=(I,\omega_I)$ in $E(A)$. Then there exists a surjection $\alpha:P\surj I$ such that $\omega_I$ is induced from 
$(\alpha, \chi)$.
\item $P$ has a unimodular element if and only if $e(P,\chi)=0$ in $E(A)$.

\item The canonical map $E(A)\longrightarrow E(\mathcal{A})$ is injective.
\end{enumerate}
\et

\rmk\label{orien} Let $R$ be a ring containing $\mathbb{Q}$ of
dimension $n$ with $\text{ht}(\mathcal{J}(R))\geq 2$.  Let $I\subset
A=R[T,T^{-1}]$ be an ideal of height $n\geq 3$ and let
$\omega_I:(A/I)^n\surj I/I^2$ be a local orientation of $I$. Let
$\theta\in GL_n(A/I)$ be such that $det(\theta)=\ol f$. Then
$\omega_I\circ\theta$ is another local orientation of $I$, which we denote
by $\ol f \omega_I$. On the other hand, if $\omega_I$ and $\widetilde
\omega_I$ are two local orientations of $I$, then by (\ref{bhat1}), it
is easy to see that $\widetilde \omega_I=\ol f \omega_I$ for some unit
$\ol f\in A/I$.
\medskip

The following result is due to Bhatwadekar and Raja Sridharan
\cite[Theorem 4.5]{brs2} in domain case. Same proof works in general
(See \cite[Proposition 5.3]{d2}).

\bt\label{br} Let $R$ be an affine algebra over an algebraically
closed field of characteristic $0$ with dim$(R)= n$.  Let $P$ be a
projective $R[T]$-module of rank $n$ with trivial determinant.
Suppose there is a surjection $\alpha:P\surj I$, where $I$ is an ideal
of $R[T]$ of height $n$ such that $I$ is generated by $n$
elements. Then $P$ has a unimodular element.  \et
\medskip

One can obtain the following result from \cite{su2} and \cite[Lemma 2.4]{ra2},

\bp\label{weak} Let $R$ be an affine algebra 
over an algebraically closed field of characteristic $0$ with dim$(R)=
n\geq 3$. Let $L$ be a projective $R$-module of rank one. Then the
canonical map $E(R,L)\longrightarrow E_0(R,L)$ is an isomorphism.  \ep

The following result is a consequence of local-global principle
for Euler class groups \cite[Theorem 4.17]{dz}. 

\bp\label{clg} Let $R$ be a ring containing $Q$ with $dim(R)= n \geq
3$. Let $P$ be a projective $R[T]$-module of rank $n$ with determinant
$L[T]$. Suppose that $P/TP$ as well as $P\ot R_m[T]$ have a unimodular
element for every maximal ideal $\mm$ of $R$. Then $P$ has a
unimodular element.  \ep

 \proof Let $\chi:L[T]\iso \wedge^n(P)$
be an isomorphism and $\alpha:P \surj I$ be a surjection, where
$I\subset R[T]$ is an ideal of height $n$.  Let
$e(P,\chi)=(I,\omega_I)$ in $E(R[T],L[T])$, where $(I,\omega_I)$ is
obtained from the pair $(\alpha,\chi)$.
As $P\ot R_m[T]$ has a unimodular element for every maximal ideal
$\mm$ of $R$, the image of $e(P,\chi)$ in
$\prod_{\mathfrak{m}}E(R_{\mathfrak{m}}[T],L_{\mathfrak{m}}[T])$ is
trivial. By \cite[Theorem 4.17]{dz} the following sequence of groups
$$0\longrightarrow E(R,L)\longrightarrow E(R[T],L[T])\longrightarrow
\prod_{\mathfrak{m}}E(R_{\mathfrak{m}}[T],L_{\mathfrak{m}}[T])$$ 
is
exact. Therefore, there exists $(J,\omega_J)\in E(R,L)$ such that
$(JR[T],\omega_J\ot R[T])=e(P,\chi)$ in $E(R[T],L[T])$.  Then we have
$e(P/TP,\chi\ot R[T]/(T))=(J,\omega_J)$ in $E(R,L)$. Since $P/TP$ has
a unimodular element, by \cite[Corollary 4.4]{brs1}, $(J,\omega_J)=0$
in $E(R,L)$. Consequently , $e(P,\chi)=0$ in $E(R[T],L[T])$. Hence, by
\cite[Corollary 4.15]{dz}, $P$ has a unimodular element.  \qed
\medskip

The following result is due to Swan \cite[Lemma 1.3]{sw}. We give a
proof due to Murthy.

\bl\label{swan} Let $R$ be a ring and $P$ be a projective
$R[T,T^{-1}]$-module. Let $f(T)\in R[T]$ be monic such
that $P_{f(T)}$ is extended from $R$. Then $P$ is extended
from $R[T^{-1}]$.  \el 

\proof If degree of $f(T)$ is $m$, then
$f(T)T^{-m}=1+T^{-1}g_1(T^{-1}):=g(T^{-1})$. Consider the following
Cartesian diagram.
 $$
\diagram
R[T^{-1}] \rrto\ar@{->>}[d] && R[T^{-1}]_{T^{-1}} (=R[T,T^{-1}])\ar@{->>}[d]\\
R[T^{-1}]_{g(T^{-1})}\rrto &\ar[r] &  R[T,T^{-1}]_{g(T^{-1})}=R[T,T^{-1}]_{f(T)}
\enddiagram 
$$

Since $P_{f(T)}$ is extended from $R$, $P_{g(T^{-1})}$ is extended
from $R[T^{-1}]_{g(T^{-1})}$. Using a standard patching argument, there
exists a projective $R[T^{-1}]$-module $Q$ such that $Q\ot_{R[T^{-1}]}
R[T,T^{-1}]\simeq P$.  This completes the proof.  \qed

\section{Main results}

In this section, we shall prove our main results stated in the
introduction. 

\bt \label{MI} Let $R$ be a ring of dimension $n$ containing an
infinite field. Let $P$ be a projective $R[T_1,\cdots,T_k]$-module of
rank $n$. Assume $P_{f(T_k)}$ has a unimodular element for some monic
polynomial $f(T_k)$ in the variable $T_k$. Then $P$ has a unimodular
element.  \et

\proof When $k=1$, we are done by \cite[Theorem 3.4]{brs2}. Assume
$k\geq 2$ and use induction on $k$. Since $(P/T_1P)_f$ has a
unimodular element, by induction on $k$, $P/T_1P$ has a unimodular
element. Also $(P\otimes R(T_1)[T_2,\ldots,T_k])_f$ has a unimodular
element and dim $(R(T_1))=n$. So again by induction on $k$, $P\otimes
R(T_1)[T_2,\ldots,T_k]$ has a unimodular element. Since $P$ is
finitely generated, there exists a monic polynomial $g\in R[T_1]$ such
that $P_g$ has a unimodular element. Applying (\ref{uni}), we get that
$P$ has a unimodular element.
\qed
\medskip

The next result is due to Das-Sridharan \cite[Theorem 3.4]{drs} when $L=R$.

\bp\label{one var} Let $R$ be a ring of {\bf even}
dimension $n$ containing $\mathbb{Q}$. Let $P$ be a projective
$R[T]$-module of rank $n$ with determinant $L$.  Suppose there are
surjections $\alpha:P\surj I$ and $\phi :L\op R[T]^{n-1} \surj I$,
where $I$ is an ideal of $R[T]$ of height $n$.  Assume further that
$P/TP$ has a unimodular element. Then $P$ has a unimodular element.
\ep

\proof
Let $R_{red}=R/{\mathbf{n}(R)}$, where $\mathbf{n}(R)$ is the nil
radical of $R$.  It is easy to derive that $P$ has a unimodular
element if and only if $P\ot R_{red}$ has a unimodular
element. Therefore, without loss of generality we may assume that $R$
is reduced. We divide the proof into two steps. \smallskip

\noindent
{\bf Step 1:} Assume that $L$ is extended from $R$.  Since $P/TP$ has
a unimodular element, in view of (\ref{clg}), it is enough to prove
that if $\mm$ is a maximal ideal of $R$, then $P\ot R_{\mm}[T]$ has a
unimodular element.  Note that $P\ot R_{\mm}[T]$ has trivial
determinant and $\alpha\otimes R_{\mm}[T]:P\otimes R_{\mm}[T]\surj
IR_{\mm}[T]$ is a surjection. Using surjection $\phi\otimes
R_{\mm}[T]$, we get that $IR_{\mm}[T]$ is generated by $n$ elements.
Applying \cite[Theorem 3.4]{drs} to the ring $R_{\mm}[T]$, it follows
that $P\ot R_{\mm}[T]$ has a unimodular element. Hence we are done.
\smallskip

\noindent
{\bf Step 2:} 
Assume that $L$ is not necessarily extended from $R$. Since $R$ is reduced, we can find a ring $R\inj S\inj Q(R)$ such that 
\begin{enumerate}
\item
The projective $S[T]$-module $L\ot_{R[T]} S[T]$ is extended from $S$,
\item $S$ is a finite $R$-module,
\item the canonical map $\text{Spec}(S)\hookrightarrow \text{Spec} (R)$
is bijective, and 
\item
for every $\mathfrak{P}\in \text{Spec}(S)$, the inclusion $R/(\mathfrak{P}\cap R)\longrightarrow
S/\mathfrak{P}$ is birational.
\end{enumerate}
Since $L\ot_{R[T]} S[T]$ is extended from $S$, by Step 1, $P\ot_{R[T]} S[T]$ has a unimodular element.
Applying \cite[Lemma 3.2]{b1}, we conclude that $P$ has a unimodular element. 
\qed

\medskip

The following result extends (\ref{one var}) to polynomial ring
$R[T_1,\cdots,T_k]$.

\bt\label{many var}
Let $R$ be a ring  of {\bf even} dimension $n$ containing $\mathbb{Q}$. Let $P$ be a projective $A=R[T_1,\cdots,T_k]$-module of rank $n$
 with determinant $L$.  Suppose there are surjection $\alpha:P\surj I$ 
and $\phi : L\op A^{n-1} \surj I$, where $I\subset A$ is an
ideal of height $n$. Assume further that $P/{(T_1,\cdots,T_k)P}$ 
has a unimodular element. Then $P$ has a unimodular element.
\et

\proof When $k=1$, we are done by (\ref{one var}). Assume $k\geq2$ and
use induction on $k$.  We use ``bar'' when we move modulo $(T_k)$.
Since $R$ contains $\mathbb{Q}$, replacing $X_k$ by $X_k-\gl$ for some 
$\gl\in \mathbb Q$, we
may assume that ht $(\ol I)\geq n$.

Consider the surjection $\ol \alpha:\ol P\surj \ol I$ induced from
$\alpha$.  If $\text{ht}(\ol I)>n$, then $\ol I=R[T_1,\cdots,T_{k-1}]$
and hence $\ol P$ contains a unimodular element. Assume $\text{ht}(\ol
I)=n$. Since $\ol P/(T_1,\cdots,T_{k-1})\ol P= P/{(T_1,\cdots,T_k)P}$
has a unimodular element and $\phi$ induces a surjection $\ol \phi:
\ol{L\op A^{n-1}} \surj \ol I$, by induction hypothesis, $\ol P$ has a
unimodular element.  By similar arguments, we get that $P/T_{k-1}P$
has a unimodular element.

Write $\mathcal{A}= R(T_k)[T_1,\cdots,T_{k-1}]$ and $P\ot
\mathcal{A}=\widehat{P}$. We have surjections $\alpha\ot \mathcal{A}:
\widehat P\surj I\mathcal{A}$ and $\phi \otimes\mathcal{A} :(L\ot
\mathcal{A})\op \mathcal{A}^{n-1} \surj I\mathcal{A}$. If
$I\mathcal{A}=\mathcal{A}$, then $\widehat P$ has a unimodular
element.  Assume $\text{ht}(I\mathcal{A})=n$.  Since
$\widehat{P}/T_{k-1}\widehat{P} = P/T_{k-1}P \otimes \mathcal A$, we
get that $\widehat{P}/T_{k-1}\widehat{P}$ and hence $\widehat{P}
/(T_1,\cdots,T_{k-1})\widehat{P}$ has a unimodular element.  By
induction on $k$, $\widehat P$ has a unimodular element. So,
there exists a monic polynomial $f\in R[T_k]$ such that $P_f$ has a
unimodular element.  Applying (\ref{uni}), we get that $P$
has a unimodular element.   \qed
\medskip

Now we prove (\ref{laurent0}) mentioned in the introduction.

\bt\label{laurent} Let $R$ be a ring containing $\mathbb{Q}$ of {\bf
  even} dimension $n$ with $\text{ht}(\mathcal{J}(R))\geq 2$.  Let $P$
be a projective $R[T,T^{-1}]$-module of rank $n$ with trivial
determinant. Assume there exists a surjection $\alpha:P \surj I$, where
$I\subset R[T,T^{-1}]$ is an ideal of height $n$ such that $I$ is
generated by $n$ elements. Then $P$ has a unimodular element.  \et

\proof Let $\chi:R[T,T^{-1}]\iso \wedge^n(P)$ be an isomorphism. Let
$e(P,\chi)=(I,\omega_I)\in E(R[T,T^{-1}])$ be obtained from the pair
$(\alpha,\chi)$.  By (\ref{global} (3)), it is enough to prove that
$e(P,\chi)=(I,w_I)=0$ in $E(R[T,T^{-1}])$. 

 Suppose $\omega_I$ is given by $I=(g_1,\cdots, g_n)+I^2$. 
Also $I$ is generated by $n$ elements, say, $I=(f_1,\cdots,f_n)$. By (\ref{orien}), there exists
$\tau\in GL_n(R[T,T^{-1}]/I)$ such that $(\ol f_1,\cdots,\ol f_n)=(\ol g_1,\cdots,\ol g_n)\tau$, where ``bar'' denotes
reduction modulo $I$. 
Let $U,V\in R[T,T^{-1}]$ be such that $det(\tau)=\ol U$ and $UV=1$ modulo $I$. 
\smallskip

\noindent
{\bf Claim:} The unimodular row $(U,f_1,\cdots,f_n)$ over
$R[T,T^{-1}])$ is completable.
\smallskip

First we show that the theorem follows from the claim.
 Since $n$ is even, by (\ref{gene}), there exists $\eta\in M_n(R[T,T^{-1}])$ with $det(\eta)=\ol V$ such that if 
  $(f_1,\cdots,f_n)\eta=(h_1,\cdots,h_n)$, then $I=(h_1,\cdots,h_n)$.
 Further, $(\ol g_1,\cdots,\ol g_n)\tau\circ \ol \eta=(\ol h_1,\cdots,\ol h_n).$
Note $\tau\circ \ol \eta\in SL_n(R[T,T^{-1}]/I)$. 

Let $\mathcal{A}$ be the ring obtained from $R[T,T^{-1}]$ by inverting
all special monic polynomials.  Since
dim$(\mathcal{A}/{I\mathcal{A}})=0$, we have
$SL_n(\mathcal{A}/{I\mathcal{A}})=E_n(\mathcal{A}/{I\mathcal{A}})$.
Let $\Theta\in SL_n(\mathcal{A})$ be a lift of $(\tau\circ\ol
\eta)\otimes \mathcal A$.  Then $I\mathcal{A}=(h_1,\cdots,h_n)
\Theta^{-1}$ and $(\ol h_1,\cdots,\ol h_n)\ol {\Theta^{-1}}=(\ol
g_1,\cdots,\ol g_n)$. In other words, $(I\mathcal{A}, \omega_I\ot
\mathcal{A})=0$ in $E(\mathcal{A})$.  Since
$\text{ht}(\mathcal{J}(R))\geq 2$, by (\ref{global} (4)), the
canonical map $E(R[T,T^{-1}])\longrightarrow E(\mathcal{A})$ is
injective.  Therefore, $(I,\omega_I)=0$ in $E(R[T,T^{-1}])$. This
completes the proof. So we just need to prove the claim.
\medskip

\noindent {\bf Proof of the claim:} Let $Q$ be the stably free
$R[T,T^{-1}]$-module associated to the unimodular row
$(U,f_1,\cdots,f_n)$.  If $S$ is the set of all monic polynomials in
$R[T]$, then $S^{-1}R[T,T^{-1}]=R(T)$. Applying (\ref{stably}),
we have $P_S$ is free and hence there exists a monic polynomial $f\in
R[T]$ such that $P_f$ is free. In particular, $P_f$ is extended from
$R$.  Hence by (\ref{swan}), there exists a projective
$R[T^{-1}]$-module $Q_1$ such that $Q\simeq Q_1\ot_{R[T^{-1}]}
R[T,T^{-1}]$.  Since $(Q_1\op R[T^{-1}])_{T^{-1}} \simeq Q\op R[T,T^{-1}]\simeq R[T,T^{-1}]^{n+1}$,
by Quillen and Suslin
\cite{q,su}, we get $Q_1\op R[T^{-1}]\simeq R[T^{-1}]^{n+1}$.
By (\ref{lg}), $Q_1$ is extended from $R$, say, $Q_1\simeq Q_2\ot_R R[T^{-1}]$. Since $\text{ht}(\mathcal{J}(R))\geq 1$,  
 by Bass \cite{bass},  $Q_2$ is free. 
Hence $Q$ is free, i.e. $(V,f_1,\cdots,f_n)$ is completable. 
This proves the claim.
\qed
\medskip

Next result is the converse of (\ref{laurent}).

\bt
Let $R$ be a ring containing $\mathbb{Q}$ of {\bf even} dimension $n$ with $\text{ht}(\mathcal{J}(R)\geq 2$.
Let $P$ be a projective $R[T,T^{-1}]$-module of rank $n$ with trivial determinant. Let $I\subset R[T,T^{-1}]$ be an ideal of height $n$
 which is generated by $n$ elements.  Suppose that $P$ has a unimodular element. Then there exists a surjection $\alpha:P\surj I$.
\et
\proof
Let $\chi:R[T,T^{-1}]\iso \wedge^n(P)$ be an isomorphism
and $e(P,\chi)\in E(R[T,T^{-1}])$ be obtained from the pair $(P,\chi)$.
 Since $P$ has a unimodular element, by (\ref{global} (3)), 
$e(P,\chi)=0$ in $E(R[T,T^{-1}])$. 

Let $I=(f_1,\cdots,f_n)$ and $\omega_I$ be the orientation of $I$
induced by $f_1,\cdots,f_n$.  Then $(I,\omega_I)=0=e(P,\chi)$ in
$E(R[T,T^{-1}])$.  By (\ref{global} (2)), there exists a surjection
$\alpha:P\surj I$ such that $(I,\omega_I)$ is induced from
$(\alpha,\chi)$.  \qed
\medskip

\bt\label{discrete} Let $R$ be a ring of {\bf even} dimension $n$
containing $\mathbb{Q}$ and $D=R[X,Y]/(XY)$. Let $P$ be a projective
$D$-module of rank $n$ with determinant $L$.  Suppose there are
surjections $\alpha:P\surj I$ and $\phi :L\oplus D^{n-1}\surj I$,
where $I$ is an ideal of $D$ of height $n$. Assume further that
$P/(X,Y)P$ has a unimodular element. Then $P$ has a unimodular
element.  \et

\proof Without loss of generality, we may
assume that $R$ is reduced. We continue to denote the
images of $X$ and $Y$ in $D$ by $X$ and $Y$. We give proof in two steps.

\noindent
{\bf Step 1:} Assume $L$ is extended from $R$.  Let ``bar'' and
``tilde'' denote reductions modulo $(Y)$ and $(X)$ respectively.  We
get surjections $\ol \alpha:\ol P\surj \ol I$ and
$\widetilde\alpha:\widetilde P\surj \widetilde I$ induced from
$\alpha$. Note that $\text{ht}(\ol I)\geq n$ and $\text{ht}(\widetilde
I)\geq n$.  If $\text{ht}(\ol I)> n$, then $\ol I=R[X]$ and $\ol P$
has a unimodular element. Assume $\text{ht}(\ol I)=n$. Since $\ol
D=R[X]$ and $\ol P/X\ol P=P/(X,Y)P$ has a unimodular element, by
(\ref{one var}), $\ol P$ has a unimodular element. Similarly,
$\widetilde P$ has a unimodular element.
 
Let $p_1\in Um(\ol P)$. Then $\widetilde p_1\in Um(\widetilde{\ol P})$. 
 Since $det(\widetilde P)=\widetilde L$ is extended from $R$, it follows from \cite[Theorem 5.2, Remark 5.3]{blr} that
 the natural  map $Um(\widetilde P)\twoheadrightarrow Um(\widetilde{\ol P})$ is surjective. 
 Therefore there exists $p_2\in Um(\widetilde P)$ such that $\ol p_2=\widetilde p_1$. 
Since the following square of rings 
 $$
\diagram
D \rrto\ar@{->>}[d] && R[X]\ar@{->>}[d]\\
R[Y]\rrto &\ar[r] &  R
\enddiagram 
$$
is Cartesian with vertical maps surjective, $p_1\in Um(\ol P)$ and 
$p_2\in Um(\widetilde P)$
will patch up to give a unimodular element of $P$. 
\smallskip

\noindent
{\bf Step 2:}
Assume $L$ is not necessarily extended from $R$. Since $R$ is reduced, by \cite[Lemma 3.7]{kz2}, 
there exists a ring $S$ such that
\begin{enumerate}
 \item $R\hookrightarrow S\hookrightarrow Q(R)$,
\item $S$ is a finite $R$-module,
\item $R\hookrightarrow S$ is subintegral and
\item $L\ot_R S $ is extended from $S$.
\end{enumerate}
Note that $R\inj S$ is actually a finite subintegral extension. Hence,
by \cite[Lemma 2.11]{kz2}, $R[X,Y]/(XY)\inj S[X,Y]/(XY)$ is also
subintegral. Since $L\ot_R S$ is extended from $S$, by Step 1, 
$P\ot S[X,Y]/(XY)$ has a unimodular element. Since $R$ contains
$\mathbb{Q}$, by \cite[Corollary 3.14]{dz1}, $P$ has
a unimodular element.  \qed
\medskip

Using \cite[Theorem 5.2, Remark 5.3]{blr} and following the proof of
(\ref{discrete}), we can prove the following result.

\bt Let $R$ be a ring of {\bf even} dimension $n$ containing
$\mathbb{Q}$ and $D=R[Y_1,\ldots,Y_m]/Y_1(Y_2,\ldots,Y_m)$. Let $P$
be a projective $D$-module of rank $n$ with determinant $L$.  Suppose
there are surjections $\alpha:P\surj I$ and $\phi :L\oplus
D^{n-1}\surj I$, where $I$ is an ideal of $D$ of height $n$. Assume
further that $P/(Y_1,\ldots,Y_m)P$ has a unimodular element. Then $P$ has a
unimodular element.  \et


\section{Applications}

In this section, we give some applications of results proved
earlier. When $L=R$, the following result is proved in \cite[Theorem
  4.5]{brs2}.

\bt \label{A1} Let $R$ be an affine algebra over an algebraically
closed field of characteristic $0$ with dim$(R)=n\geq 3$.  Let $P$ be
a projective $R[T]$-module of rank $n$ with determinant $L$.  Suppose
$I$ is an ideal of $R[T]$ of height $n$ such that there are two
surjections $\alpha:P\surj I$ and $\phi:L\oplus R[T]^{n-1} \surj I$.
Then $P$ has a unimodular element.  \et

\proof Let ``bar'' denote reduction modulo $(T)$.  Then $\alpha$
induces a surjection $\ol \alpha: \ol P \surj I(0)$.  We can assume
ht$(I(0))\geq n$.  If $\text{ht}(I(0))>n$, then $I(0)=R$ and $\ol P$
has a unimodular element.  Assume $\text{ht}(I(0))=n$ and fix an
isomorphism $\chi:L\iso \wedge^n(P)$. Let $e(\ol P,\ol
\chi)=(I(0),\omega_{I(0)})$ in $E(R,\ol L)$ be induced from $(\ol
\alpha,\ol\chi)$.  Since $\ol\phi:\ol L\oplus R^{n-1}\surj I(0)$ is a
surjection, by (\ref{weak}), $e(\ol P,\ol\chi)=0$ in $E(R,\ol L)$. By
\cite[Corollary 4.4]{brs1}, $P/TP$ has a unimodular element. Rest of
the proof is exactly as in (\ref{one var}) using (\ref{br}).  \qed
\medskip

\bc\label{discrete1} Let $R$ be an affine algebra over an
algebraically closed field of characteristic $0$ with dim$(R)=n\geq 3$
and $D=R[X,Y_1\ldots,Y_m]/X(Y_1,\ldots,Y_m)$.  Let $P$ be a projective
$D$-module of rank $n$ with determinant $L$. Suppose $I$ is an ideal
of $D$ of height $n$ such that there are two surjections
$\alpha:P\surj I$ and $\phi:L\oplus D^{n-1} \surj I$.  Then $P$ has
a unimodular element.  \ec

\proof Follow the proof of (\ref{discrete}) and use (\ref{A1}). \qed
\medskip

\bt \label{A2} Let $R$ be an affine algebra over an algebraically
closed field of characteristic $0$ such that dim$(R)=n\geq 4$ is {\bf
  even}.  Let $P$ be a projective $A=R[T_1,\cdots,T_k]$-module of rank
$n$ with determinant $L$. Suppose $I$ is an ideal of $A$ of height $n$
such that there are two surjections $\alpha:P\surj I$ and
$\phi:L\oplus A^{n-1} \surj I$.  Then $P$ has a unimodular
element. \et

\proof By (\ref{weak}) and \cite[Corollary 4.4]{brs1}, we get
$P/(T_1,\cdots,T_{k})P$ has a unimodular element. By
(\ref{many var}), we are done.  \qed
\medskip

In the following special case, we can remove the restriction on
dimension in (\ref{many var}).  This result improves \cite[Theorem
  4.4]{brs2}.
 
\bp Let $R$ be a ring of dimension $n\geq 2$ containing $\mathbb{Q}$
and $P$ be a projective $A=R[T_1,\cdots,T_k]$-module of rank $n$ with
determinant $L$. Suppose $\mathcal{M}\subset A$ is a maximal ideal of
height $n$ such that there are two surjections
$\alpha:P\surj \mathcal M$ and $\phi:L\oplus A^{n-1} \surj \mathcal M$.
Then $P$ has a unimodular element.  \ep

\proof Assume $k=1$. If $\mathcal{M}$ contains a monic polynomial
$f\in R[T_1]$, then $P_f$ has a unimodular element. Using
\cite[Theorem 3.4]{brs2}, $P$ has a unimodular element.  Assume
$\mathcal{M}$ does not contain a monic polynomial.  Since $T_1\notin
\mathcal{M}$, ideal $\mathcal{M}+ (T_1)=R[T_1]$. By surjection
$\alpha$, we get $P/T_1P$ has a unimodular element. Follow the
proof of (\ref{one var}) and use \cite[Theorem 4.4]{brs2} to conclude
that $P$ has a unimodular element. We are done when $k=1$.

Assume $k\geq 2$ and use
induction on $k$.
For rest of the proof, follow the proof of (\ref{many var}). 
\qed


\medskip



\begin{thebibliography}{}
\bibitem[B]{bass} H. Bass, K-theory and stable algebra, {\it
  Publ. Math. Inst. Hautes \'{E}tudes Sci.} {\bf 22} (1964), 5-60.

\bibitem[Bh 1]{b1} S. M. Bhatwadekar, Inversion of monic polynomials and existence of unimodular elements (II), 
{\it Math. Z.} {\bf 200} (1989), 233-238.

\bibitem[Bh 2]{b2} S. M. Bhatwadekar, Cancellation theorems for projective
  modules over a two dimensional ring and its polynomial extensions,
  {\it Compositio Math.} {\bf 128} (2001), 339-359.

\bibitem[B-R 1]{br1} S. M. Bhatwadekar, Amit Roy, Stability theorems for
  overrings of polynomial rings, {\it Invent.Math} {\bf 68} (1982),
  117-127.

\bibitem[B-R 2]{br2} S.M. Bhatwadekar and A. Roy, Some theorems about
  projective modules over polynomial rings. {\it J. Algebra} {\bf 86} (1984),
150-158.


\bibitem[B-RS 1]{brs1} S. M. Bhatwadekar and  Raja Sridharan,
Euler class group of a Noetherian ring, {\it Compositio Math.} 122 (2000), 183-222.

\bibitem[B-RS 2]{brs2} S. M. Bhatwadekar and  Raja Sridharan, 
On a question of Roitman, {\it  J. Ramanujan Math. Soc.} {\bf 16} (2001), 45-61.


\bibitem[B-L-Ra]{blr} S. M. Bhatwadekar, H. Lindel and Ravi A. Rao,
  The Bass-Murthy question: Serre dimension of Laurent polynomial
  extensions, {\it Invent. Math.} {\bf 81} (1985), no. 1, 189-203.

\bibitem[D]{d2} M. K. Das, The Euler class group of a polynomial algebra II, {\it J. Algebra} {\bf 299} (2006), 94-114.

\bibitem[D-RS]{drs} M. K. Das, Raja Sridharan, The Euler class groups of polynomial rings and unimodular elements in projective modules, 
{\it J. Pure Applied Algebra} {\bf 185} (2003), 73-86.

\bibitem[D-Z 1]{dz1} M. K. Das and Md. Ali Zinna, On invariance of the
  Euler class group under a subintegral base change, {\it J. Algebra}
  {\bf 398} (2014), 131-155.

\bibitem[D-Z 2]{dz} M. K. Das and Md. Ali Zinna, The Euler class group
  of a polynomial algebra with coefficients in a line bundle, {\it
    Math. Z.} {\bf 276} (2014), 757-783.


\bibitem[K]{k} M. K. Keshari, Euler class group of a Laurent
  polynomial ring: Local case, {\it J. Algebra} {\bf 308} (2007),
  666-685.

\bibitem[K-Z 1]{kz} M. K. Keshari and Md. Ali Zinna, Unimodular
  elements in projective modules and an analogue of a result of
  Mandal, {\it J. Commutative Algebra} (to appear).

\bibitem[K-Z 2]{kz2} M. K. Keshari and Md. Ali Zinna, Efficient
  generation of ideals in a discrete Hodge algebra, {\it J. Pure
    Applied Algebra} (to appear).



\bibitem[Ma]{ma} Satya Mandal, Basic Elements and Cancellation over
  Laurent Polynomial Rings, {\it J. Algebra} {\bf 79} (1982), 251-257.


\bibitem[P]{p} B. Plumstead, The conjecture of Eisenbud and Evans,
  {\it Amer. J. Math} {\bf 105} (1983), 1417-1433.

\bibitem[Q]{q} D. Quillen, Projective modules over polynomial rings,
  {\it Invent. Math.} {\bf 36} (1976), 167-171.

\bibitem[R]{rao} Ravi A. Rao, The Bass-Quillen conjecture in dimension
  $3$ but characteristic $\neq 2,3$ via a question of A. Suslin, {\it
  Invent. Math.} {\bf 93} (1988), 609-618.

\bibitem[Ra 1]{ra} Raja Sridharan, Non-vanishing sections of algebraic
  vector bundles, {\it J. Algebra} {\bf 176} (1995), 947-958.

\bibitem[Ra 2]{ra2} Raja Sridharan, Projective modules and complete
  intersections, {\it K-Theory} {\bf 13} (1998), 269-278.


\bibitem[Se]{se} J.P. Serre, Sur les modules projectifs, {\it
  Semin. Dubreil-Pisot} {\bf 14} (1960-61).

\bibitem[Su 1]{su} A. A. Suslin, Projective modules over polynomial
  ring are free, {\it Soviet Math. Dokl.} {\bf 17} (1976), 1160-1164.

\bibitem[Su 2]{su2} A.A. Suslin, A cancellation theorem for projective modules over affine algebras, {\it Sov. Math. Dokl.} {\bf 18} 
(1977), 1281-1284.

\bibitem[Sw]{sw} R. G. Swan, Projective Modules over Laurent Polynomial Rings, {\it Transactions of the AMS} {\bf Vol. 237} (1978), 111-120.

\bibitem[Wi]{wi} A. Wiemers, Some Properties of Projective Modules over Discrete Hodge Algebras, {\it J. Algebra} {\bf 150} (1992), 402-426.

\end{thebibliography}
\end{document}